\documentclass[12pt]{article} 
 \usepackage{dsfont}
 \usepackage{amsmath}
 \usepackage{amsthm}
 \usepackage{txfonts,bm}
 \usepackage{latexsym}
 \usepackage[colorlinks,citecolor=blue,urlcolor=blue,linkcolor=blue]{hyperref}

 \newtheorem{theo}{\small\bf Theorem}[section]
 \newtheorem{lem}{\small\bf Lemma}[section]
 
 \newtheorem{rem}{\small\bf Remark}[section]
 \newenvironment{REM}{\begin{rem} \rm}{\end{rem}}
 \newtheorem{defi}{\small\bf Definition}[section]
 
 \newtheorem{cor}{\small\bf Corollary}[section]

 \numberwithin{equation}{section}

 \newcommand{\be}{\begin{equation}}
 \newcommand{\ee}{\end{equation}}
 \newcommand{\E}{\operatorname{\mathds{E}}}
 
 \newcommand{\Var}{\operatorname{Var}}

 \newcommand{\RR}{\mathds{R}}

 \newcommand{\CC}{\mathds{C}}

 \newcommand{\ud}{{\rm d}}

\makeatletter
\renewcommand\@biblabel[1]{#1.}
\makeatother

 \newcommand{\Bin}{\Big(\hspace{-.6ex}
 \begin{array}{c} n \\ k \end{array}
 \hspace{-.6ex}\Big)}

 \newcommand{\BinB}{\Big(\hspace{-.6ex}
 \begin{array}{c} n \\ k-1 \end{array}
 \hspace{-.6ex}\Big)}

 \newcommand{\BinC}{\Big(\hspace{-.6ex}
 \begin{array}{c} n+1 \\ k \end{array}
 \hspace{-.6ex}\Big)}

 \title{ \Large\bf On the inversion of the Laplace transform
 \\
 {\it In Memory of Dimitris Gatzouras}
 \vspace*{.9em}
 }

 \author{\large\bf
  Nickos Papadatos
  \vspace*{-3ex}
  }

  \date{\normalsize
  Department of Mathematics, National and Kapodistrian
  University of Athens,
  Panepistemiopolis, 157 84 Athens, Greece.
 }

\begin{document}

 \maketitle

 \thispagestyle{empty}

 \begin{abstract}
 \noindent
 The Laplace transform is a useful and powerful analytic tool
 with applications to several areas of applied mathematics,
 including differential equations, probability and statistics.
 Similarly to the inversion of the Fourier
 transform, inversion
 formulae for the Laplace transform are of central importance;
 such formulae are old and well-known
 (Fourier-Mellin or Bromwich integral,
 Post-Widder inversion).
 The present work is motivated from an elementary statistical
 problem, namely, the unbiased estimation of a parametric function
 of the scale in the basic model of a random sample from
 exponential distribution.
 The form of the uniformly
 minimum variance unbiased estimator of a parametric
 function $h(\lambda)$, as well as its variance,
 are obtained as series in Laguerre
 polynomials and the corresponding Fourier coefficients,
 and a particular application of this result yields
 a novel inversion formula for the Laplace transform.
 \end{abstract}
 {\footnotesize
 }

 \section{Introduction}
 \label{sec.1}

 For a function $u:(0,\infty)\to \RR$, its Laplace transform is
 defined by the integral
 \be
 \label{LT}
 \phi(\lambda)=\int_0^\infty \exp(-\lambda x) u(x) \ud x,
 \ee
 provided that there exists $\lambda_0\geq 0$ such that
 \[
 \int_0^\infty \exp(-\lambda_0 x) |u(x)| \ud x<\infty.
 \]

 There is a second version of the Laplace transform,
 related to probability measures $\mu$
 supported in (a subset of) $[0,\infty)$, namely,
 \be
 \label{LT_measures}
 \phi_{\mu}(\lambda)=\int_{[0,\infty)} \exp(-\lambda x) \ud \mu(x);
 \ee
 it is just a matter of notation to express $\phi_\mu(\lambda)$
 as $\E \exp(-\lambda X)$ where the nonnegative random variable $X$ has
 distribution function $F(x)=\mu([0,x])$, $x\geq 0$, and $\E$
 denotes expectation. In this setup, $\phi_{\mu}(\lambda)$
 is denoted as $M_{X}(-\lambda)$ and it is called {\it the moment
 generating function of $X$}. It is clear that formulae
 \eqref{LT} and \eqref{LT_measures} coincide if $X$ has a density $u$
 (w.r.\ to Lebesgue measure on $[0,\infty)$).
 An inversion formula for the probabilistic version (\ref{LT_measures}) can
 be found in Billingsley (1995) or Schilling {\it et al} (2012), and it is based
 on an ingenious application of the law of large numbers. The
 formula can be written as ($x>0$)
 \be
 \label{probabilistic.inversion}
 \mu\Big([0,x)\Big)+\frac{1}{2}
 \mu\Big(\{x\}\Big)=\lim_N \sum_{k=0}^{N} \frac{(-1)^k}{k!}
 \Big(\frac{N}{x}\Big)^k \phi_\mu^{(k)}\Big(\frac{N}{x}\Big).
 \ee

 Regarding \eqref{LT}, it is known from Lerch (1903) that the transformation
 $u\rightarrow \phi$ is one to one.
 Furthermore, there are two well-known inversion formulae for \eqref{LT},
 namely, the Fourier-Mellin or Bromwich integral
 (see Boas (1983), Cohen (2007)),
 \be
 \label{FMB}
 u(x)=\frac{1}{2\pi i} \lim_{T\to\infty}
 \int_{\gamma-i T}^{\gamma+iT} \exp(sx)\phi(s) \ud s,
 \ee
 where $\gamma\geq 0$ is greater than the real part of every pole
 of (the analytic extension of) $\phi$, and the Post (1930)
 or
 Post-Widder formula
 (see Widder (1946), Post (1930), Cohen (2007)),
 \be
 \label{PW}
 u(x)=\lim_n
 \frac{(-1)^n}{n!}\left(\frac{n}{x}\right)^{n+1}\phi^{(n)}\Big(\frac{n}{x}\Big).
 \ee
 The above inversions hold under some mild restrictions, e.g.,
 \eqref{FMB} is satisfied for almost all $x\in(0,\infty)$ (clearly, this is the best
 we can expect, but the formula in itself is complicated and, so,
 inconvenient
 for purposes  of computation, as can be seen when applied to trivial exemplary cases),
 and \eqref{PW} holds at continuity points of $u$, provided
 that $u$ is smooth in pieces and that the growth of $|u|$
 at infinity is at most of exponential order.

 The present work is motivated from an elementary statistical
 inference
 problem which, at a first glance, seems to be unrelated to Laplace inversion.
 The problem is to find the minimum variance unbiased
 estimator of a given parametric function $h(\lambda)$,
 based on a random sample $X_1,\ldots,X_n$
 from Exp$(\lambda)$, with $\lambda>0$ unknown, or, more generally, from
 $\Gamma(a,\lambda)$
 with $a>0$ fixed and known and $\lambda>0$ an unknown
 parameter
 (for the definitions see Section \ref{sec.2}).
 The main result provides necessary and sufficient conditions
 on $h$ so that a solution of this problem exists, and shows
 that the solution (whenever exists)
 can be presented as a series in Laguerre
 polynomials,
 \be
 \label{Laguerre}
 L_n(x)=\sum_{k=0}^n (-1)^k
 \Big(\hspace{-.6ex}
 \begin{array}{c} n \\ k \end{array}
 \hspace{-.6ex}\Big)
 \frac{x^k}{k!}.
 \ee
 A particular application of the main result yields
 a novel inversion formula for the Laplace transform; see
 Section \ref{sec.3}.

 \section{On the best unbiased estimator of a parametric function
 of the scale parameter in exponential/gamma models}
 \label{sec.2}

 \subsection{Preliminaries and a simple parametric inference problem}
 \label{ssec.2.1}
 The probability density of the exponential distribution, Exp($\lambda$), is given by
 \[
 f_{\lambda}(x)=\lambda\exp(-\lambda x), \ \ x>0,
 \]
 while the Gamma distribution, $\Gamma(a,\lambda)$, has
 probability density
 \be
 \label{f lambda}
 f_{\lambda}(x)=\frac{\lambda^a}{\Gamma(a)}x^{a-1}\exp(-\lambda x), \ \ x>0,
 \ee
 where $a>0$ and $\lambda>0$ are positive constants,
 so that Exp$(\lambda)\equiv\Gamma(1,\lambda)$.

 From now on, we suppose that $a>0$ is known (given), and we
 assume that $\lambda>0$ is the (unique) unknown parameter
 to be estimated from the data. More generally, we wish
 to estimate an arbitrary parametric function $h(\lambda)$ by
 using a suitable choice of an estimator
 \[
 T=T(X_1,\ldots,X_n),
 \]
 where
 $T$ is a real valued measurable function with domain $(0,\infty)^n$
 and $X_1,\ldots, X_n$ are iid (independent, identically distributed)
 random variables with density \eqref{f lambda}.
 Of course, the actual value of $T$ (when $X_1=x_1,\ldots,X_n=x_n$)
 must not vary with $\lambda$, but $T$ may depends on $n$ or $a$
 (since both are fixed and known).

 So, the problem can be formulated as follows:
 \vspace*{1em}

 \noindent
 {\small\bf Problem 1.} Let $h(\lambda):(0,\infty)\to\RR$ be a given
 (arbitrary)
 parametric function and suppose that $X_1,\ldots,X_n$ are iid with density
 \eqref{f lambda}. Under what conditions on $h$ is
 it possible to find an estimator $T=T(X_1,\ldots,X_n)$
 such that
 \begin{enumerate}
 \item $\E_{\lambda} T=h(\lambda)$ \ \ for all $\lambda>0$,
 and
 \item $\E_{\lambda} T^2 < \infty$ \ \ for all
 $\lambda>0$?
 \end{enumerate}
 And, in case that such a $T$ exists,
 how can we obtain the best possible estimator for $h$?
 \vspace{1em}

 An estimator satisfying condition 1 is called {\it unbiased};
 as we shall see, unbiasedness restricts the class of
 possible estimators in such a way that the family of
 permitted  parametric functions $h$ is quite narrow.
 Condition 2 means that
 $T\in \cap_{\lambda>0}L^2(\mu_n(\lambda))$, where $\mu_n(\lambda)$
 is the product probability measure of $(X_1,\ldots,X_n)$
 on $[0,\infty)^{n}$. Then, provided $\E_{\lambda}T=h(\lambda)$,
 the quantity $\E_{\lambda} (T-h(\lambda))^2$
 can be written as $\Var_{\lambda} T = \E_\lambda T^2-h(\lambda)^2$, and it is called
 the {\it variance} of the estimator $T$. Even if $T$ is not
 unbiased, the quantity $\E_\lambda (T-h(\lambda))^2$ is called MSE
 (mean squared error), and it is the most important measure
 of closeness between a point estimator $T(X_1,\ldots,X_n)$
 and a parametric function
 $h(\lambda)$, traditionally used in statistics for a long time.
 The subscript $\lambda$ in $\E_\lambda$ and $\Var_\lambda$
 denotes that the true probability measure of the $X_i$'s
 is as in \eqref{f lambda}.

 It is clear that, if we restrict ourselves to the class of unbiased estimators,
 those
 with smaller variance are preferable. In the plausible case
 that we can pick an
 estimator $T^{*}$ satisfying
 \begin{enumerate}
 \item $\E_{\lambda} T^*=h(\lambda)$ \ \ for all $\lambda>0$,
 \item $\Var_{\lambda} T^* < \infty$ \ \ for all
 $\lambda>0$, and
 \item for any unbiased estimator $T$ and for all $\lambda>0$,
 $\Var_{\lambda} T^*\leq \Var_{\lambda} T$,
 \end{enumerate}
 it follows that this is the best we can do. Such an estimator
 $T^*$ is then
 called {\it uniformly minimum variance unbiased estimator}
 (UMVUE for short), and this is what we could name as {\it best}.
 In order to be able to obtain the UMVUE it is necessary and
 sufficient that the class
 \[
  {\cal T}_{h}=
  \{T: T \ \mbox{is an unbiased estimator for $h(\lambda)$
  with finite variance (for all $\lambda>0$)} \}
 \]
 is nonempty.
 This follows from one of the most fundamental result in parametric inference,
 adapted to the present particular case of Gamma distributions
 with $a$ known. Indeed, the following is true; see Lehmann and Gasella
 (1998).

 \begin{theo} {\rm (Rao-Blackwell / Lehmann-Scheff\'{e}).} \
 \label{theo.RB LS}
 Let $X_1,\ldots,X_n$ be a random sample from \eqref{f lambda}
 with $\lambda>0$ unknown and $a>0$ known. Let also
 $h:(0,\infty)\to\RR$ be a parametric function, and suppose that
 ${\cal T}_h$ is nonempty. Set
 $X=X_1+\cdots+X_n$.
 Then,
 \medskip

 \noindent
 {\rm (i)} The conditional probability distribution of $(X_1,\ldots,X_n)$
 given $X$ is independent of $\lambda$.
 \medskip

 \noindent
 {\rm (ii)} For any $T\in{\cal T}_h$, the
 {\rm (unique w.p.\ 1)} UMVUE is given by the conditional expectation
 \[
 T^*(X):=\E\Big(T(X_1,\ldots,X_n) \ \Big| \ X\Big).
 \]

 \noindent
 {\rm (iii)} Equivalently, the
 UMVUE  of $h(\lambda)$ is the unique {\rm (w.p.\ 1)} unbiased
 estimator in ${\cal T}_{h}$ which is a function of $X$,
 $u=u(X)$. Hence,
 $u(X)=\E\Big(T(X_1,\ldots,X_n) \ \Big| \ X\Big)=T^*(X)$.
 \end{theo}

 \begin{REM}
 \label{rem.substitution}
 It is well-known that the distribution of $X=X_1+\cdots+X_n$ is
 $\Gamma(n a,\lambda)$. In view of Theorem \ref{theo.RB LS},
 and substituting $a$ for $n a>0$,
 Problem 1 reduces to an equivalent, much simper,
 reformulation, as follows.
 \end{REM}

 \noindent
 {\small\bf Problem 2.} Let $h(\lambda):(0,\infty)\to\RR$ be a given
 (arbitrary)
 parametric function and suppose that $X$ is a random variable
 having probability density \eqref{f lambda}, with $a>0$ fixed and known, and $\lambda>0$
 an unknown parameter.
 Under what conditions on $h$ does
 the UMVUE $u=u(X)$
 of $h(\lambda)$ exists for all $\lambda$?
 And, in case that it exists,
 how can we obtain its form?
 \vspace{1em}

 Since, by definition, $\E_{\lambda} \psi(X)=\int_0^{\infty} f_{\lambda}(x)\psi(x) \ud x $
 for arbitrary measurable $\psi$,
 the imposed condition of a finite second moment on $u$
 {\it for all
 $\lambda$} implies that
 \be
 \label{L_2 lambda}
 \int_{0}^\infty x^{a-1} \exp(-\lambda x) u(x)^2 \ud x<\infty.
 \ee
 In other words,
 $u\in L^2(\lambda)$ for all $\lambda>0$,
 where $L^2(\lambda)$ is the Lebesgue space of functions
 $u:(0,\infty)\to\RR$
 satisfying \eqref{L_2 lambda}.
 Thus, it is reasonable to define
 \be
  \mbox{$L_0^2:=\bigcap_{\lambda>0} L^2(\lambda)$}.
 \ee

 We can rewrite the unbiasedness restriction $\E_{\lambda} u(X)=h(\lambda)$
 as
 \be
 \label{analytic}
 \frac{\Gamma(a) h(\lambda)}{\lambda^a} =\int_0^\infty
 x^{a-1} \exp(-\lambda x) u(x) \ud x, \ \ \lambda>0.
 \ee
 It is then obvious that the rhs of \eqref{analytic} defines a holomorphic
 function in the right half-plane $\CC^+=\{\lambda\in \CC:
 \Re(\lambda)>0\}$ whenever $u\in L^2_0$. This means that
 the function
 $\lambda\rightarrow \Gamma(a) \lambda^{-a}
 h(\lambda)$ is holomorphic, and hence, $h(\lambda)$
 must be holomorphic in $\CC^+$. This already imposes
 a serious restriction to the allowable parametric functions,
 e.g., it is necessary that $h\in C^{\infty}(0,\infty)$; in fact,
 the analytic extension of $h$ should have no singularities
 in the right half-plane. As a simple example, for the $C^{\infty}(0,\infty)$
 parametric function
 $h(\lambda)=1/(\lambda^2-2\lambda+2)$,
 no unbiased estimator exist (for all $\lambda>0$),
 because of the poles $1\pm i$ of $h$.
 However, regarding Problem 2,  the analyticity of $h$ is
 not sufficient to provide a positive answer. To see this, it suffices
 to observe that for $u\in L^2_0$,
 $\int_0^{\infty} x^{a-1} \exp(-\lambda x)  u(x)\ud x\to 0$
 as $\lambda\to+\infty$, by dominated convergence. Then, any
 holomorphic function $h$ that growths faster than $\lambda^a$ at
 infinity, e.g.\ $h(\lambda)=\exp(\lambda)$,
 cannot be written as the expectation of some $u\in L^2_0$; see \eqref{analytic}.

 \subsection{Results}
 \label{ssec.results}
 We are now in a position to state the main results.
 \begin{theo}
 \label{theo.main}
 Assume that $X$ is a random variable
 with probability density $f_{\lambda}$
 as in \eqref{f lambda}, with $\lambda>0$ unknown.
 For a given parametric function $h(\lambda)$, its UMVUE $u(X)$
 exists in $L^2_0$ if and only if
 the following two conditions are satisfied.
 \\
 {\rm (1)} The function $h$ can be extended to an holomorphic function
 in $\CC^+$, and \\
 {\rm (2)} For any $\lambda>0$,
 \be
 \label{Fourier el-2}
 \sum_{n=0}^\infty \beta_n(\lambda)^2<\infty,
 \ee
 where
 \be
 \label{Fourier_constants}
 \beta_n(\lambda)=\frac{(-1)^n}{\sqrt{n! [a]_n}}
 \Bigg(\lambda\frac{\ud^n}{\ud\lambda^n}\Big[\lambda^{n-1}
 h(\lambda)\Big]\Bigg);
 \ee
 here, $[a]_n=\prod_{j=0}^{n-1} (a+j)=\Gamma(a+n)/\Gamma(a)$
 denotes the ascending factorial {\rm (Pochhammer symbol)}.
 \end{theo}

 \begin{theo}
 \label{theo.main.constants}
 Let $h$ be a parametric function
 satisfying {\rm (1)} 
 of Theorem {\rm \ref{theo.main}}.
 For fixed $\lambda>0$ define the
 function
 \be
 \label{f_simplified}
 H_\lambda(y):=h\Big(\frac{\lambda}{1-y}\Big), \ \ |y|<1.
 \ee
 Then, an alternative simplified form of the constants $\beta_n(\lambda)$
 in {\it \eqref{Fourier_constants}} is given by
 \be
 \label{Fourier_constants.simple}
 \beta_n(\lambda)=
 \frac{(-1)^n H_\lambda^{(n)}(0)}{\sqrt{n! [a]_n}}.
 \ee
 \end{theo}

 \begin{theo}
 \label{theo.u}
 Assume that {\rm (1)} and {\rm (2)} of Theorem {\rm
 \ref{theo.main}}
 are satisfied and fix $\lambda_0>0$.
 Then, the function $u(x)$ for which $u(X)$
 is the UMVUE of $h(\lambda)$ is given by
 \be
 \label{u}
 u(x)=\sum_{n=0}^\infty \beta_n(\lambda_0) q_{n;\lambda_0}(x)
 \ee
  where $\{q_{n;\lambda_0} (x)\}_{n=0}^{\infty}$ is the complete orthonormal polynomial
  system corresponding to the weight function
  $f_{\lambda_0}$,
  with the convention that each $q_{n;\lambda_0}$ is of degree $n$
  and with strictly positive leading coefficient. The series converges
  a.e.\ on $(0,\infty)$ and in $L^2 (\lambda)$ for every $\lambda\geq\lambda_0$,
  and the resulting function $u(x)$, given by \eqref{u}, is independent of
  the choice of $\lambda_0$. Furthermore, for any $\lambda>0$,
  the variance of the UMVUE is given
  by
  \be
  \label{var UMVUE}
  \Var_{\lambda} u(X)=\sum_{n=1}^{\infty} \beta_n(\lambda)^2,
  \ee
  where the constants $\beta_n(\lambda)$ are completely determined from $h(\lambda)$;
  see \eqref{Fourier_constants} or
  \eqref{Fourier_constants.simple}.
 \end{theo}

 \noindent
 {\bf Example 1.}
 We compare the expression \eqref{var UMVUE} with the
 classical information inequality, namely, the famous Cram\'{e}r-Rao
 (CR)
 lower bound ($\mbox{LB}_{CR}$). Since, as is well-known, the regularity conditions are satisfied
 for $f_{\lambda}$, the bound states that
 for a random sample $X_1,\ldots,X_n$ (of size $n$) from $f_\lambda$, and for
 any unbiased estimator $T=T(X_1,\ldots,X_n)$ of $h(\lambda)$,
 the inequality
 $\Var_{\lambda} T\geq h'(\lambda)^2/(nI(\lambda)):=\mbox{LB}_{CR}$
 is satisfied; here, $I(\lambda)$ is the Fisher information,
 defined as
 \[
 I(\lambda):=\E_{\lambda}\Bigg[\Bigg(\frac{\partial}{\partial
 \lambda} \log f_{\lambda}(X_1)\Bigg)^2\Bigg]=\frac{a}{\lambda^2}.
 \]
 Thus, the CR-bound reads as $\Var_{\lambda} T\geq \lambda^2
 h'(\lambda)^2/(na)$. On the other hand, the series
 expansion \eqref{var UMVUE} (applied with $na$ in place of $a$; see
 Problems 1 and 2) yields
 \[
 \Var_{\lambda} u(X)=\sum_{m=1}^{\infty} \frac{\lambda^2}{m! [n a]_m}
 \Bigg(\frac{\ud^m }{\ud
 \lambda^m}\Big[\lambda^{m-1}h(\lambda)\Big]\Bigg)^2.
 \]
 Since $u(X)$ is the UMVUE and thus,
 $\Var_{\lambda} T\geq \Var_{\lambda}u(X)$
 for any unbiased estimator $T$,
 it is clear that the CR-bound is implied by the preceding
 series, on just keeping its first term.

 \subsection{Proofs}
 \label{ssec.proofs}
 We first state some auxiliary lemmas.
 \begin{lem}
 For $x>0$, $a>0$ and $\lambda>0$,
 \label{lem.key}
 \be
 \label{key}
 \frac{\ud ^n}{\ud x^n}
 \Big[ x^n f_{\lambda}(x)\Big]=\lambda \frac{\ud^n}{\ud\lambda^n}
 \Big[\lambda^{n-1}f_\lambda(x)\Big], \ \ n=0,1,2,\ldots \ .
 \ee
 \end{lem}
 \noindent
 {\small\bf Proof.} By Leibnitz formula and \eqref{f lambda} it is easily seen that
 both sides of \eqref{key}
 are equal to
 \[
 \Gamma(a+n) f_{\lambda}(x)  \sum_{k=0}^n (-1)^k \Bin
 \frac{(\lambda x)^k}{\Gamma(a+k)}.
 \]

 \begin{lem}
 {\rm (Rodrigues' formula; see Afendras and Papadatos (2015)).} \
 For $x>0$, $a>0$ and $\lambda>0$,
 \label{lem.Rodriguez}
 \be
 \label{Rodriguez}
 \frac{\ud ^n}{\ud x^n} \Big[ x^n f_{\lambda}(x)\Big]=
 (-1)^n  \sqrt{[a]_n n!} f_\lambda(x)q_{n;\lambda}(x),
 \ \ n=0,1,2,\ldots,
 \ee
 where $\{q_{n;\lambda}(x)\}_{n=0}^{\infty}$ is the complete
 orthonormal system with respect to $f_{\lambda}$,
 standardized so that $q_{n;\lambda}$
 has degree $n$ and positive leading coefficient. The polynomials
 $q_{n;\lambda}$ satisfy the orthogonality condition
 \vspace*{-2ex}
 \[
 \E_{\lambda} \Big[q_{n;\lambda}(X) q_{m;\lambda}(X)\Big]
 =\mbox{$\int_0^{\infty} f_{\lambda}(x) q_{n;\lambda}(x)
 q_{m;\lambda}(x) \ \ud x$}
 =\left\{
 \begin{array}{ccc}
 1 & \mbox{if} & n=m;
 \\
 0 & \mbox{if} & n\neq m.
 \end{array}
 \right.
 \]
 \end{lem}

 One important observation is that, as
 \eqref{Rodriguez} and \eqref{key} show, we may
 produce
 the orthonormal set $q_{n;\lambda}$ by
 differentiate w.r.\ to the parameter $\lambda$,
 instead of $x$.; more precisely,
 \be
 \label{orthonormal}
 q_{n;\lambda}(x)=\frac{(-1)^n }{\sqrt{n! [a]_n}f_{\lambda}(x)}
 \frac{\ud^n}{\ud x^n}
 \Big[ x^n f_{\lambda}(x) \Big]=
 \frac{(-1)^n  }{\sqrt{n! [a]_n}f_{\lambda}(x)}
 \Bigg(\lambda\frac{\ud^{n}}{\ud {\lambda^{n}}}
 \Big[\lambda^{n-1}f_{\lambda}(x)
 \Big]\Bigg).
 \ee
 Thus, \eqref{orthonormal} obtains the following
 \begin{cor}
 \label{cor.OPS}
 For $x>0$, $a>0$, $\lambda>0$ and $n\in\{0,1,\ldots\}$,
  \be
 \label{orthogonal.lambda}
 q_{n;\lambda}(x) f_{\lambda}(x)=
 \frac{(-1)^n  }{\sqrt{n! [a]_n}}
 \Bigg(\lambda\frac{d^{n}}{d {\lambda^{n}}} \Big[\lambda^{n-1}f_{\lambda}(x)
 \Big]\Bigg).
 \ee
 \end{cor}

 \noindent
 We now proceed to verify the results of Theorems
 \ref{theo.main}--\ref{theo.u}.

 Assume first that the UMVUE of $h(\lambda)$ is $u(X)$, and suppose that it has
 finite variance
 for all $\lambda>0$.
 Multiplying the equation $\E_{\lambda} u(X)=h(\lambda)$ by $\lambda^{n-1}$
 and then taking $n$ derivatives w.t.\ to $\lambda$, we subsequently obtain
 \begin{eqnarray}
 \nonumber
 h(\lambda) & = & \int_{0}^\infty f_{\lambda}(x) u(x) \ud x,
 \\
 \nonumber
 \lambda^{n-1} h(\lambda) & = & \int_{0}^\infty \lambda^{n-1} f_{\lambda}(x) u(x) \ud x,
 \\
 \nonumber
 \frac{\ud^{n}}{\ud\lambda^n}
 \Big[\lambda^{n-1} h(\lambda)\Big] & = &
 \int_{0}^\infty \left(\frac{\ud^n}{\ud\lambda^n}
 \Big[\lambda^{n-1} f_\lambda(x)\Big]\right) u(x) \ud x,
 \\
 \nonumber
 \lambda \frac{\ud^{n}}{\ud\lambda^n}
 \Big[\lambda^{n-1} h(\lambda)\big] & = &
 \int_{0}^\infty \left(
 \lambda \frac{\ud^n}{\ud\lambda^n}\Big[\lambda^{n-1} f_\lambda(x)\Big]\right) u(x) \ud x,
 \\
 \label{coefficients.derivatives}
 \frac{(-1)^n}{\sqrt{n! [a]_n}}
 \Bigg(\lambda\frac{\ud^{n}}{\ud\lambda^n}
 \Big[\lambda^{n-1} h(\lambda)\big]\Bigg) & = &
 \int_{0}^\infty q_{n;\lambda}(x) f_{\lambda}(x) u(x) \ud x;
 \end{eqnarray}
 note that the differentiation can be passed under the integral
 sign, due to the assumed (squared) integrability of $u$ with
 respect to $f_\lambda$ for all $\lambda>0$. We conclude from
 \eqref{coefficients.derivatives} that the constants $\beta_n(\lambda)$
 of \eqref{Fourier_constants} are the Fourier coefficients of
 $u$ with respect to the orthonormal polynomial system
 $\{q_{n;\lambda}\}_{n=0}^{\infty}$, corresponding to the weight function
 $f_\lambda$. It should be noticed that the orthonormal polynomial system
 corresponding to a probability measure (having finite moments of any order)
 is unique, apart from a possible multiplication of each polynomial
 by $\pm 1$. Moreover, since our system $\{q_{n;\lambda}\}_{n= 0}^{\infty}$
 is complete in $L^2(\lambda)$, see Afendras {\it et al} (2011),
 Parseval's identity yields
 \[
  \E_{\lambda} u(X)^2=
  \int_{0}^{\infty} f_{\lambda}(x) u(x)^2 \ud x =
  \sum_{n=0}^\infty
  \beta_n(\lambda)^2<\infty.
 \]
 Thus, assuming that $u\in L^2_0$, and since
 $\beta_0(\lambda)=\E_{\lambda}u(X)=h(\lambda)$,
 we get
 \[
 \Var_{\lambda} u(X) = \sum_{n=1}^\infty
  \beta_n(\lambda)^2, \ \ \mbox{ for all $\lambda>0$}.
 \]

 Conversely, assume that $h$ is holomorphic in $\CC^{+}$ and that the series
 in \eqref{Fourier_constants} is finite for all $\lambda>0$.
 Then we may define the function $u(x;\lambda)$ as
 \be
 \label{u(x;lambda)}
 u(x;\lambda):=\sum_{n=0}^\infty \beta_n(\lambda) q_{n;\lambda}(x),
 \ \ x>0,
 \ee
 where, by Riesz-Fisher, the series converges in $L^2(\lambda)$, that is,
 \[
  \int_0^\infty \Big(u_N(x;\lambda)-u(x;\lambda)\Big)^2
  f_{\lambda}(x) \ud x \to 0, \ \ N\to\infty,
 \]
 with $u_N(x;\lambda)=\sum_{n=0}^N \beta_n(\lambda)
 q_{n;\lambda}(x)$. It remains to show that
 the limiting function $u(x;\lambda)$ does not depend on
 $\lambda$,
 and that it is the unique UMVUE of $h(\lambda)$.
 To this end, choose a fixed $\lambda_0>0$ with $\lambda_0<\lambda$ and write
 \be
 \label{L2_limit lambda_0}
 u_N(x)=\sum_{n=0}^N \beta_n(\lambda_0)
 q_{n;\lambda_0}(x), \ \ x>0, \ N=0,1,2,\ldots \ .
 \ee
 From the convergence of the series \eqref{Fourier el-2}
 (with $\lambda=\lambda_0$)
 it is easily seen that $u_N(x)$ is Cauchy $L^2(\lambda_0)$, and hence, it converges
 (in the norm of $L^2(\lambda_0)$) to a function $u(x)\in L^2(\lambda_0)$.
 Moreover, is easy to check that for any $\lambda\geq \lambda_0$, we can
 find a constant $C_{\lambda}=C(\lambda,\lambda_0)$ such that
 $f_{\lambda}(x)\leq C_\lambda f_{\lambda_0}(x)$ for all $x>0$. This
 implies that $u_N$ is also Cauchy $L^2(\lambda)$ for any fixed $\lambda\geq \lambda_0$;
 indeed,
 if $\epsilon>0$ is arbitrary,
 we can find $N(\epsilon)$
 such that $\int_0^{\infty} \left(u_{N_1}(x)-u_{N_2}(x)\right)^2 f_{\lambda_0}(x) \ud
 x< \epsilon/C_{\lambda}$ for all $N_1,N_2>N(\epsilon)$ and, then,
 \[
 \int_0^{\infty} \left(u_{N_1}(x)-u_{N_2}(x)\right)^2 f_\lambda(x) \ud
 x
 \leq C_{\lambda}\int_0^{\infty} \left(u_{N_1}(x)-u_{N_2}(x)\right)^2 f_{\lambda_0}(x) \ud
 x<\epsilon.
 \]
 The preceding argument verifies that the limiting function $u$,
 defined as the $L^2(\lambda_0)$-limit of the sequence in \eqref{L2_limit lambda_0},
 belongs to $L^2(\lambda)$ for all $\lambda\geq \lambda_0$, in symbols,
 $u(x)\in \bigcap_{\lambda\geq\lambda_0}L^2(\lambda)$. From
 the orthogonality of the polynomials
 $q_{n;\lambda_0}$ ($n\geq 1$) and $q_{0;\lambda_0}\equiv 1$
 we immediately see that
 $\E_{\lambda_0} u_N(X)=\beta_0(\lambda_0)=h(\lambda_0)$, and
 clearly, this is also true for $u$, i.e.,
 $\E_{\lambda_0}u(X)=h(\lambda_0)$.
 However, the situation is different
 when $\lambda\neq \lambda_0$, that is,
 the expectation of $u_N(X)$ w.r.\ to $f_\lambda$ varies with
 both $N$ and $\lambda$.
 More precisely, since $q_{0;\lambda}(x)\equiv 1$,
 \[
 \E_{\lambda} u_N(X)
 =
 h(\lambda_0)+
 \sum_{n=1}^N \beta_n(\lambda_0)
 \E_{\lambda} q_{n;\lambda_0}(X)
 , \ \  \ N=1,2,\ldots \ \ \lambda>0.
 \]
 On the other hand, we have shown that for $\lambda\geq \lambda_0$,
 $\E_{\lambda}\Big(u_N(X)-u(X)\Big)^2\to 0$,
  so that, by the Cauchy-Schwarz inequality,
  \[
  \Big|\E_{\lambda}u_N(X)-\E_{\lambda}u(X) \Big|
  \leq \E_{\lambda} \big|u_N(X)-u(X)\big| \leq
  \Big(\E_{\lambda} \big|u_N(X)-u(X)\big|^2\Big)^{1/2}\to 0.
  \]
 It follows that $\E_\lambda u(X)=\lim_N \E_{\lambda} u_N(X)$.
 Hence, the expectation of $u(X)$ w.r.\ to $f_\lambda$
 can be obtained
 as the limit of the expectations of $u_N(X)$ (w.r.\ to
 $f_{\lambda}$). Next, we see that the calculation of
 $\E_{\lambda}u_N(X)$
 requires
 evaluation of the expectations $\E_{\lambda} q_{n;\lambda_0}(X)$,
 that is, integrals of the polynomials $q_{n;\lambda_0}(x)$
 w.r.\ to a different weight function ($f_{\lambda}$ instead of $f_{\lambda_0}$),
 under which these polynomials are no longer orthogonal.

 In order to calculate $\E_{\lambda} q_{n;\lambda_0}(X)$ we proceed
 as follows. We have
 \begin{eqnarray*}
 \E_{\lambda} q_{n;\lambda_0}(X)
 &=&
 \int_0^\infty \frac{f_{\lambda}(x)}{f_{\lambda_0}(x)}
 f_{\lambda_0}(x) q_{n;\lambda_0}(x) \ud x
 \\
 &=&
 \Big(\frac{\lambda}{\lambda_0}\Big)^{a}
 \int_0^\infty
 f_{\lambda_0}(x) \exp\Big(-(\lambda-\lambda_0)x\Big)q_{n;\lambda_0}(x)
 \ud x.
 \end{eqnarray*}
 The last integral can be viewed as the $n$-th Fourier coefficient of the
 bounded $C^{\infty}(0,\infty)$ function
 $w(x):= \exp\Big(-(\lambda-\lambda_0)x\Big)$, $x>0$, with respect to
 the corresponding orthonormal polynomial system $\{q_{n;\lambda_0}\}_{n=0}^\infty$
 in $L^2(\lambda_0)$. On the other hand, it is known that the same Fourier coefficients
 can be conveniently obtained by using the identity
 (see Afendras and Papadatos (2015), Afendras {\it et al} (2011))
 \[
 \E_{\lambda_0} \Big[q_{n;\lambda_0}(X) w(X)\Big]=
 \frac{1}{\sqrt{n![a]_n}}
 \E_{\lambda_0} \Big[X^n w^{(n)}(X)\Big],
 \]
 provided $\E_{\lambda_0} \Big[X^n \Big(w^{(n)}(X)\Big)^2\Big]<\infty$.
 Since $w^{(n)}(x)= (\lambda_0-\lambda)^n
 \exp\Big(-(\lambda-\lambda_0)x\Big)$
 is a bounded function of $x$, because $\lambda\geq \lambda_0$,
 we can apply the preceding
 formulae
 to deduce
 \[
 \E_{\lambda} q_{n;\lambda_0}(X)
 =
 \Big(\frac{\lambda}{\lambda_0}\Big)^{a}
 \frac{(\lambda_0-\lambda)^n}{\sqrt{n![a]_n}}
 \E_{\lambda_0} \Big[X^n \exp\Big(-(\lambda-\lambda_0)X\Big)\Big].
 \]
 A straightforward computation now yields
 \[
 \E_{\lambda_0} \Big[X^n \exp\Big(-(\lambda-\lambda_0)X\Big)\Big]
 =
 \frac{\lambda_0^a}{\Gamma(a)}
 \int_0^{\infty} x^{n+a-1} \exp(-\lambda x)
  \ud x
    =
 [a]_n \frac{\lambda_0^a}{\lambda^{n+a}},
 \]
 and thus,
 \[
 \E_{\lambda} q_{n;\lambda_0}(X)
 =(-1)^n
 \sqrt{\frac{[a]_n}{n!}}
 \left(1-\frac{\lambda_0}{\lambda}\right)^n.
 \]
 Recalling that $\beta_n(\lambda_0)$ is given by
 \eqref{Fourier_constants} with $\lambda=\lambda_0$, we
 have
 \begin{eqnarray}
 \E_{\lambda} u_N(X)
 & = &
 \sum_{n=0}^N
 \frac{(-1)^n\lambda_0}{\sqrt{n! [a]_n}}
 \frac{\ud^n}{\ud \lambda^n}\Big[\lambda^{n-1}
 h(\lambda)\Big]\Bigg|_{\lambda=\lambda_0} \left(
 (-1)^n
 \sqrt{\frac{[a]_n}{n!}}
 \left(1-\frac{\lambda_0}{\lambda}\right)^n\right)
 \nonumber
 \\
 &=&
 \label{Taylor}
 \sum_{n=0}^N
 \frac{1}{n!}
 \left(1-\frac{\lambda_0}{\lambda}\right)^n
 \Bigg\{\lambda_0 \frac{\ud^n}{\ud \lambda^n}\Big[\lambda^{n-1}
 h(\lambda)\Big]\Bigg|_{\lambda=\lambda_0}\Bigg\}.
 \end{eqnarray}
 Though the preceding formula appears to be quite complicated at a first glance
 (e.g., it seems that it is not an easy task to obtain its limiting value as $N\to\infty$),
 this is not the case. In fact,
 \eqref{Taylor} represents a Taylor development around $y=0$ of the
 function $H_{\lambda_0}(y):=h\Big(\frac{\lambda_0}{1-y}\Big)$, $|y|<1$.
 Recall that $h(\lambda)$ has been assumed to be holomorphic in $\Re(\lambda)>0$,
 so that $H_{\lambda_0}(y)$ is analytic in the open disc $|y|<1$.
 Writing $H_{\lambda_0}^{(n)}(y)$ for $\frac{\ud^n}{\ud y^n}H_{\lambda_0}(y)$, we
 shall verify below the equality
 \be
 \label{Taylor H}
 H_{\lambda_0}^{(n)}(0)=\Bigg\{\lambda_0 \frac{\ud^n}{\ud \lambda^n}\Big[\lambda^{n-1}
 h(\lambda)\Big]\Bigg|_{\lambda=\lambda_0}\Bigg\}, \ \ n=0,1,\ldots \ .
 \ee
 Assuming for a while that \eqref{Taylor H} is valid,
 and substituting it to \eqref{Taylor}, we obtain
 the simple formula
 \[
  \E_{\lambda} u_N(X)=
  \sum_{n=0}^N
 \frac{H_{\lambda_0}^{(n)}(0)}{n!}
 \left(1-\frac{\lambda_0}{\lambda}\right)^n.
 \]
 Since $|1-\lambda_0/\lambda|<1$ (for $\lambda>\lambda_0/2$),
 we conclude from Taylor's theorem that
 $\E_{\lambda} u_N(X)\to H_{\lambda_0}(1-\lambda_0/\lambda)=h(\lambda)$.
 Thus, $\E_{\lambda} u(X)=\lim_N \E_{\lambda} u_N(X)=h(\lambda)$,
 and this verifies that $u(X)$ is the (unique) UMVUE of $h(\lambda)$,
 for every $\lambda\geq \lambda_0$.
 [To see uniqueness, repeat
 the previous
 construction  with $\lambda_1$ in place of $\lambda_0$.
 Then, as we showed, the produced estimator $\widetilde{u}(X)$
 will satisfy $\E_{\lambda} \widetilde{u}(X)=h(\lambda)=\E_{\lambda}u(X)$
 for all $\lambda\geq \max\{\lambda_0,\lambda_1\}$, so it must be
 identical to $u(X)$.] \ Furthermore, \eqref{Fourier_constants}
 shows that $u$ has the same Fourier coefficients
 as the function $u(x;\lambda)$ defined by \eqref{u(x;lambda)};
  thus $u(x)=u(x;\lambda)$ is independent of $\lambda$, and Parseval's
 identity yields \eqref{var UMVUE}. The orthogonal polynomials
 for the weight function $f_{\lambda}$ are called generalized
 Laguerre (Laguerre when $a=1$).
 The a.e.\ convergence of the Laguerre series expansion of a function
 $u\in L^2(\lambda)$ is a well-known (Carleson-Hunt type) result that can be found
 in Mackenhoupt (1970); see also Uspensky (1927) and Stempak (2000).

 It remains to show \eqref{Taylor H}. Using Leibnitz formula
 we first calculate
 \be
 \label{L.1}
 \lambda \frac{\ud ^n}{\ud \lambda^n} \Big[\lambda^{n-1} h(\lambda)\Big]
 = (n-1)! \sum_{k=1}^n \Bin  \frac{\lambda^{k}
 h^{(k)}(\lambda)}{(k-1)!}, \ \ n=1,2,\ldots,
 \ee
 while the lhs equals to $h(\lambda)$ for $n=0$.
 Next, we define $H_{\lambda}(y)=h\Big(\lambda/(1-y)\Big)$, $|y|<1$,
 so that $H_{\lambda}^{(0)}(y)=H_{\lambda}(y)$ and $H_{\lambda}^{(0)}(0)=h(\lambda)$.
 For $n=1$,
 $H_{\lambda}'(y)=\lambda
 h'\Big(\lambda/(1-y)\Big)\Big/(1-y)^2$, and
 $H_{\lambda}'(0)=\lambda h'(\lambda)$ equals to the sum in the rhs of \eqref{L.1}
 (with $n=1$). We shall prove, using induction on $n$, the formula
 (valid for $\lambda>0$, $|y|<1$)
 \be
 \label{L.2}
 H_{\lambda}^{(n)}(y) = (n-1)! \sum_{k=1}^n \Bin  \frac{\lambda^{k}
 h^{(k)}(\lambda/(1-y))}{(k-1)!(1-y)^{n+ k}}, \ \ n=1,2,\ldots,
 \ee
 which, setting $y=0$, yields the rhs of \eqref{L.1}; then,
 the substitution $\lambda\rightarrow\lambda_0$ verifies \eqref{Taylor H}.
 Noting that \eqref{L.2} is true for $n=1$, we assume that it is true
 for some $n$. Then,
 \begin{eqnarray*}
 H_{\lambda}^{(n+1)}(y)
 &=&
 \frac{\ud}{\ud y} \Bigg\{(n-1)! \sum_{k=1}^n \Bin  \frac{\lambda^{k}
 h^{(k)}(\lambda/(1-y))}{(k-1)!(1-y)^{n+k}}\Bigg\}
 \\
 &=&
  (n-1)! \sum_{k=1}^n \Bin
  \frac{\lambda^k}{(k-1)!}
  \frac{\ud}{\ud y}\Bigg\{\frac{ h^{(k)}(\lambda/(1-y))}{(1-y)^{n+ k}}\Bigg\}
 \\
 &=&
  (n-1)! \sum_{k=1}^n \Bin
  \frac{\lambda^k}{(k-1)!}
  \frac{ h^{(k+1)}(\lambda/(1-y))}{(1-y)^{n+ k}}\frac{\lambda}{(1-y)^2}
  \\
  &&
  \ \ \ \ + \ (n-1)! \sum_{k=1}^n \Bin
  \frac{\lambda^k}{(k-1)!}
  \frac{ h^{(k)}(\lambda/(1-y))}{(1-y)^{n+ k+1}} (n+k)
  \\
  &=&
  (n-1)! \sum_{k=2}^{n+1} (k-1)\BinB
  \frac{\lambda^{k}
  h^{(k)}(\lambda/(1-y))}{(k-1)!(1-y)^{n+1+ k}}
  \\
  &&
  \ \ \ \ + \ (n-1)! \sum_{k=1}^n (n+k) \Bin
  \frac{\lambda^{k}
  h^{(k)}(\lambda/(1-y))}{(k-1)!(1-y)^{n+1+ k}}
  \\
  &=&
  (n-1)! \sum_{k=1}^{n+1} \Bigg\{(k-1)\BinB
   +(n+k) \Bin
  \Bigg\}
  \frac{\lambda^{k}
  h^{(k)}(\lambda/(1-y))}{(k-1)!(1-y)^{n+1+ k}},
  \end{eqnarray*}
  where the last equality follows from $\Bin=0$
  for $k=n+1$ and $(k-1)\BinB=0$ for $k=1$.
  It is now obvious that
  \begin{eqnarray*}
  (k-1)\BinB
   +(n+k) \Bin
   &= &
   \frac{(k-1) n!}{(k-1)!(n-k+1)!}+ \frac{(n+k) n!}{k!(n-k)!}
   \\
   &=&
   \Big\{
    k(k-1)+(n+k)(n+1-k)
   \Big\}
   \frac{n!}{k!(n+1-k)!}
   \\
   &=&
   n(n+1) \frac{n!}{k!(n+1-k)!}
   \\
   &=&
   n \BinC.
  \end{eqnarray*}
 This shows that \eqref{L.2} holds with $n+1$ in place of $n$,
 and concludes the inductional argument.

 \section{A novel inversion formula of the Laplace transform}
 \label{sec.3}

 The results of Section \ref{sec.2} apply to the particular case
 where $a=1$, i.e., when $X$ follows the exponential distribution
 with parameter $\lambda>0$,  Exp$(\lambda)$, with probability density
 \be
 \label{f lambda exponential}
 f_\lambda(x)=\lambda \exp(-\lambda x), \ \ x>0.
 \ee
 In this case, Lemma \ref{lem.Rodriguez} produces the corresponding
 orthonormal
 polynomial system, namely,
 \[
 q_{n;\lambda}(x)=\sum_{k=0}^n (-1)^{n-k} \Bin \frac{(\lambda
 x)^k}{k!}.
 \]
 The preceding polynomials are functions of $\lambda x$ (this is also
 true for $a\neq 1$, since it is easily seen that
 $q_{n;\lambda}(x)=q_{n;1}(\lambda x)$).
 Hence, it is convenient to define $p_n(x)=q_{n;1}(x)$, so that
 $q_{n;\lambda}(x)=p_n(\lambda x)$.
 Then, the polynomial system $\{p_n(x)\}_{n=0}^\infty$
 is the complete orthonormal system corresponding to $f_1$ (i.e., Exp$(1)$),
 that is,
 \[
 \E \Big[p_{n}(X)p_m(X)\Big] =\int_0^{\infty} e^{-x} p_{n}(x)p_m(x)
 \ud x =
 \left\{
 \begin{array}{ccc}
 1 & \mbox{if} & n=m,
 \\
 0 & \mbox{if} & n\neq m,
 \end{array}
 \right.
 \]
 where $\E$ stands for $\E_1$.
 Traditionally, the polynomials $L_n(x)=(-1)^n p_n(x)$ (with alternating
 leading coefficients)
 are called {\it Laguerre polynomials},
 see \eqref{Laguerre}, and they are also orthonormal w.r.\ to $f_1(x)=e^{-x}$,
 $x>0$.

 Consider now Problem 2 with $a=1$. This reduces in finding the
 function
 \[
 u\in L^2_0:=\mbox{$\bigcap_{\lambda>0}$}L^2\Big((0,\infty),e^{-\lambda x}\Big)
 \]
 for which
 \[
 \E_\lambda u(X) :=\int_{0}^\infty \lambda \exp(-\lambda x)
 u(x) \ud x  =h(\lambda), \ \ \lambda>0,
 \]
 provided that $h(\lambda)$ allows such a construction. Theorem
 \ref{theo.main} provides a necessary and sufficient condition on
 $h$, namely, $h(\lambda)$ is holomorphic for
 $\lambda\in\CC^{+}=\{\lambda\in \CC: \Re(\lambda)>0\}$ satisfying
 \be
 \label{sufficient}
 \sum_{n=0}^\infty
 \left(\frac{(-1)^n}{n!}
 \Bigg(\lambda\frac{\ud^n}{\ud\lambda^n}\Big[\lambda^{n-1}
 h(\lambda)\Big]\Bigg)\right)^2<\infty, \ \ \lambda>0.
 \ee
 In view of Theorem \ref{theo.main.constants},
 the preceding condition can be rewritten as
 \be
 \label{sufficient simple}
 \sum_{n=0}^\infty
 \left(
 \frac{(-1)^n H_\lambda^{(n)}(0)}{n!}\right)^2<\infty, \ \
 \lambda>0,
 \ee
 where $H_{\lambda}(y)=h\Big(\lambda/(1-y)\Big)$, $|y|<1$.

 It is obvious that the
 equation $\E_\lambda u(X)=h(\lambda)$ can be written in terms of the Laplace
 transform of $u$, \eqref{LT},
 as
 \[
 \lambda \phi(\lambda)=\int_{0}^\infty \lambda\exp(-\lambda x) u(x) \ud x
 =\E_{\lambda} u(X)=h(\lambda).
 \]
 Hence, given the (holomorphic in $\CC^{+}$) Laplace transform $\phi$,
 one can check the validity of either \eqref{sufficient}
 or \eqref{sufficient simple} for
 $h(\lambda):=\lambda\phi(\lambda)$,
 in order to ensure that the inverse function $u(x)$
 of $\phi(\lambda)$ exists in $L^2_0$; if this is the case, then
 Theorem \eqref{theo.u} applies and $u$ is obtained
 as a Laguerre polynomial series with constants
 derived from the derivatives of $\phi$.

 Translating Theorems \ref{theo.main}--\ref{theo.u} to
 the Laplace-transform case, we obtain the following
 \begin{theo}
 \label{theo.Laplace.Inversion}
 {\rm (A)} Assume that $\phi(\lambda)$ is an
 holomorphic function of $\lambda\in\CC^+$,
 such that
 \be
 \label{sufficient.Laplace}
 \sum_{n=0}^\infty
 \left(\frac{1}{n!}
 \Bigg(\lambda\frac{\ud^n}{\ud\lambda^n}\Big[\lambda^{n}
 \phi(\lambda)\Big]\Bigg)\right)^2<\infty, \ \ \lambda>0,
 \vspace*{-1.5ex}
 \ee
 or, equivalently,
 \vspace*{-1.5ex}
 \be
 \label{sufficient simple.Laplace}
 \sum_{n=0}^\infty
 \left(
 \frac{ \Phi_\lambda^{(n)}(0)}{n!}\right)^2<\infty, \ \
 \lambda>0,
 \vspace*{-0.5ex}
 \ee
 where
 \vspace*{-0.5ex}
 \be
 \label{Phi}
 \Phi_{\lambda}(y)=\frac{\lambda}{1-y} \ \phi\Big(\frac{\lambda}{1-y}\Big),
 \ \ |y|<1.
 \vspace*{-0.5ex}
 \ee
 Then, $\phi$ is the Laplace transform of a function $u\in L^2_0$.
 Moreover, for every fixed $\lambda_0>0$,
 the inverse Laplace transform, $u$, is given by
 \vspace*{-0.5ex}
 \be
 \label{Inverse.LT}
 u(x)=\sum_{n=0}^\infty
 \frac{\Phi_{\lambda_0}^{(n)}(0)}{n!} \ L_n(\lambda_0 x),
 \vspace*{-0.5ex}
 \ee
 where the Laguerre polynomials $L_n$ are given by \eqref{Laguerre}.
 The series converges a.e.\  and in $L^2\big(\RR_+,e^{-\lambda x}\big)$
 for every $\lambda\geq \lambda_0$, and the sum of the series does not
 depend
 on the particular choice of $\lambda_0$.

 \noindent
 {\rm (B)} If $\phi$ is the Laplace transform of a function $u\in L^2_0$
 then $\phi$ is holomorphic in $\CC^{+}$ and satisfies
 \eqref{sufficient.Laplace} {\rm (equivalently, \eqref{sufficient
 simple.Laplace})}.

 \end{theo}

 Since the choice of $\lambda_0$ does not affect the validity
 of \eqref{Inverse.LT}, we
 may set $\lambda_0=1$. Then, the function $\Phi_{\lambda}$
 in \eqref{Phi} reduces to
 $\Phi_1(y)=(1-y)^{-1}\phi\big((1-y)^{-1}\big)=\Phi(y)$, say, and
 from \eqref{Inverse.LT} we obtain
 the (Taylor-like) Laplace inversion formula
 \vspace*{-0.5ex}
 \be
 \label{Taylor-like Inversion}
 u(x)=\sum_{n=0}^\infty
 \frac{\Phi^{(n)}(0)}{n!} \ L_n(x), \ \
 \mbox{ where } \ \Phi(y)=\frac{1}{1-y}\ \phi\Big(\frac{1}{1-y}\Big),
 \vspace*{-1.5ex}
 \ee
 which is valid almost everywhere in $(0,\infty)$.

 At this point we note that all inversion formulae of $\phi(\lambda)$
 provide approximating
 functions for $u(x)$ in some sense. For instance,
 \eqref{Taylor-like
 Inversion} says that
 \vspace*{-0.5ex}
 \be
 \label{Taylor_like Limit}
 u_N(x) := \sum_{n=0}^N
 \frac{\Phi^{(n)}(0)}{n!} \ L_n(x) \to u(x), \ \ \mbox{a.e.},
 \vspace*{-1ex}
 \ee
 while
 \vspace*{-0.5ex}
 \eqref{FMB} can be written in our case as
 \[
 w_N(x):=\frac{1}{2\pi i}
 \int_{1-i N}^{1+i N} \exp(sx)\phi(s) \ud s \to u(x), \ \
 \mbox{a.e.},
 \vspace*{-1ex}
 \]
 and \eqref{PW} reads as
 \vspace*{-0.5ex}
 \[
 v_N(x):=
 \frac{(-1)^N}{N!}\left(\frac{N}{x}\right)^{N+1}\phi^{(N)}\Big(\frac{N}{x}\Big)
 \to u(x) \ \ \mbox{at continuity points $x$ of $u(x)$.}
 \vspace*{-.5ex}
 \]
 Hence, it would be desirable to compare the degree of approximation of the
 preceding formulae; however, this is beyond the scope of the
 present work. We merely point out a possible advantage of the new
 inversion formula: The approximating functions $u_N$ in \eqref{Taylor_like Limit}
 are polynomials, and
 the formula becomes exact for any polynomial $u$ when $N\geq \deg(u)$.
 \bigskip

 \noindent
 {\bf Acknowledgement.} I would like to cordially thank A.\
 Giannopoulos for organizing the meeting in memory of our friend
 and colleague Dimitris, and also for providing me
 bibliographic material regarding the a.e.\
 convergence of Laguerre series.

 {
 
 }
 \end{document}